\documentclass{elsart}
\usepackage{natbib}
\usepackage{amsthm, amsmath, amsfonts}

\begin{document}

\begin{frontmatter}

\title{Structured matrices in the application of bivariate
interpolation to curve implicitization}

\author{Ana Marco \thanksref{EM1}}
\author{, Jos\'e-Javier Mart{\'\i}nez \thanksref{EM2}}

\address{Departamento de Matem\'aticas, Universidad
de Alcal\'a,}

\address{Campus Universitario, 28871-Alcal\'a de Henares (Madrid), Spain}

\thanks[EM1]{Corresponding author. E-mail: ana.marco@uah.es}

\thanks[EM2]{E-mail: jjavier.martinez@uah.es}

\begin{abstract}

A nonstandard application of bivariate polynomial interpolation is
discussed: the implicitization of a rational algebraic curve given
by its parametric equations. Three different approaches using the
same interpolation space are considered, and their respective
computational complexities are analyzed. Although the techniques
employed are usually associated to numerical analysis, in this
case all the computations are carried out using exact rational
arithmetic. The power of the Kronecker product of matrices in this
application is stressed.

\end{abstract}


\end{frontmatter}

\section{Introduction}

{\sl Curve implicitization}, which consists of finding the
implicit equation of a curve $C$ given by a rational
parametrization, is an important problem in computer aided
geometric design, and several theoretical results which help to
its solution have been developed in the fields of classical
algebraic geometry and computer algebra. Some of the methods for
the effective computation of the implicit equation are based on
interpolation.

We will consider three types of interpolation problems which share
the same interpolation space. The different ways of choosing the
interpolation nodes will lead to linear systems with very
different coefficient matrices: an unstructured one, the transpose
of a Vandermonde matrix, and the Kronecker product of two
Vandermonde matrices.

The Kronecker product structure of the matrix in the third method
will make much less expensive the computational cost of the
process [12, 22] and will introduce a high degree of {\sl
parallelism}.

In this application of interpolation {\sl exact arithmetic} is
used, so stability is not as important as efficiency and
 we do not have to worry about matters which are very
important in numerical analysis such as the ill-conditioning of
Vandermonde matrices [9]. In fact, the algorithms to be described
use techniques from numerical linear algebra but all the
computations are carried out in exact arithmetic, which allows
them to be applied with no difficulty even when polynomials of
very high degree are involved.

Finally, the function to be interpolated is itself a bivariate
polynomial, so there will be no interpolation error, and
consequently in this case the use of interpolation will not be
related to approximation theory.

\medskip

In order to make our exposition as clear and complete as possible,
we will begin with a very small example which will be useful to
illustrate the different approaches to the problem. Let us
consider the hyperbola with parametric equations
$$
(x(t), y(t)) = \Big( \frac{1+t}{2+t}, \frac{3+t}{4+t} \Big).
$$
Its implicit equation is
$$
2-3y-x+2xy=0.
$$
In the following sections it will be shown how this implicit
equation can be obtained.

The first problem to be addressed is the choice of an appropriate
interpolation space. Let us observe that the implicit equation of
any degree $n$ polynomial rational parametric curve, with the same
denominator in both components of the parametrization, is a degree
$n$ algebraic curve (see Section 15.4 of [20]). However, this is
not necessarily true when two different denominators are
considered in the parametrization. In our example, both parametric
equations have degree $1$ while the implicit equation has total
degree $2$. This is  the reason why, when considering the general
situation of two different denominators, it is more natural to
consider the {\it coordinate degree} (i.e. the degree in $x$, the
degree in $y$) instead of the total degree. The following result
(see [15] and [21]) gives us the precise form of the implicit
equation.

\medskip

{\bf Theorem 1.} {\it Let $P = \Big( x(t)=
{{u_1(t)}\over{v_1(t)}},y(t)={{u_2(t)}\over{v_2(t)}} \Big)$ be a
proper rational parametrization of the irreducible curve $C$
defined by $F(x,y)$, and let \break $gcd(u_1,v_1) = gcd(u_2,v_2) =
1$. Then $max\{deg_t(u_1),deg_t(v_1)\} = deg_y(F)$ and
 $max\{deg_t(u_2),deg_t(v_2)\} = deg_x(F)$.}

\medskip

Theorem $1$ tells us that the polynomial $F(x,y)$ defining the
implicit equation of the curve $C$ belongs to the polynomial space
$\Pi_{m,n}(x,y)$, where $m = max \{deg_t(u_2),deg_t(v_2)\}$ and $n
= max \{deg_t(u_1),deg_t(v_1)\}$. The dimension of
$\Pi_{m,n}(x,y)$ is $N=(m+1)(n+1)$, and a basis is given by
$$
\{x^i y^j | i = 0, \cdots, m; j = 0, \cdots, n\}= $$ $$ = \{1, y,
\cdots, y^n, x, xy, \cdots, xy^n, \cdots, x^m, x^my, \cdots,
x^my^n\}.
$$

Moreover $deg_x(F(x,y)) = m$ and $deg_y(F(x,y)) = n$, and
therefore there is no interpolation space $\Pi_{r,s}(x,y)$ with
$r<m$ or $s<n$ such that $F(x,y)$ belongs to $\Pi_{r,s}(x,y)$.

In addition, the selection of the interpolation space
$\Pi_{m,n}(x,y)$ is also suitable because in practice the implicit
representation of a rational parametric curve is a dense
polynomial [13].

So in the application of interpolation we are considering the
interpolation space is given and we have to choose appropriate
Lagrange interpolation nodes. As we will see, for obtaining a
linear system with a structured coefficient matrix we will need
the introduction of {\sl resultants}, a tool which is widely used
in computer algebra and has a variety of applications [4].

In the next three sections we will discuss three different
approaches to the implicitization problem by using interpolation,
and in Section 5 we will analyze their respective computational
complexities.

\bigskip
\bigskip

\section{An interpolation problem with an unstructured coefficient matrix}

The first approach, included in [11] along with other approaches
to implicitization, leads to a non-structured interpolation
problem whose solution is not unique.

Since the dimension of the linear space of solutions will be seen
to be almost surely $1$, the implicit equation can be computed in
the following way. As we know the interpolation space (it is
$\Pi_{m,n}(x,y)$), we may evaluate the parametrization $P=(x(t),
y(t))$ at some finite set of values of $t$ (for example $t=0,1,2,
\ldots$) for obtaining $N$ distinct rational interpolation nodes
$$
\{ (x_i, y_i)=(x(t_i), y(t_i)): \quad i=1, \ldots, N \}.
$$

Then we can formulate an interpolation problem using those
interpolation nodes, in this case with {\it all the interpolation
data equal to zero}. So, we have an interpolation problem in which
the corresponding linear system is a non-structured homogeneous
one, and we know from Section 1 that a nontrivial solution always
exists.

If we denote by $A$ the coefficient matrix of this homogeneous
linear system and $r = rank(A)$ then the set of all solutions of
$Ax=0$ is the nullspace of $A$, and its dimension $N-r$ is great
or equal than $1$.

Let us recall that B\'ezout's theorem states that {\it two plane
curves of degree $m+n$ without common components have at most
$(m+n)^2$ common complex points}. As a consequence of this
theorem, the probability of having $rank(A) < N-1$ is negligible,
since we are prescribing $N=(m+1)(n+1)$ interpolation conditions
at rational points. However, if that happens we can add a new
interpolation condition to find the correct implicit equation.

\medskip
In our example the rank is necessarily equal to $3$ because
$m=n=1$ and therefore $(m+n)^2=4=(m+1)(n+1)$. If we consider $t=0,
1, 2, 3$ the interpolation nodes in the corresponding order are
$$
\Big\{\Big(\frac{1}{2}, \frac{3}{4}\Big), \Big(\frac{2}{3},
\frac{4}{5}\Big), \Big(\frac{3}{4}, \frac{5}{6}\Big),
\Big(\frac{4}{5}, \frac{6}{7}\Big) \Big\},
$$
and the coefficient matrix of the homogeneous linear system $Ac=0$
is
$$
A= {\Large
\begin{pmatrix}
1 & \frac{3}{4} & \frac{1}{2} & \frac{3}{8} \\
1 & \frac{4}{5} & \frac{2}{3} & \frac{8}{15} \\
1 & \frac{5}{6} & \frac{3}{4} & \frac{5}{8} \\
1 & \frac{6}{7} & \frac{4}{5} & \frac{24}{35}
\end{pmatrix}}.
$$
The solution vector is
$$
c=(2, -3, -1, 2)^T
$$
or any multiple of it, and so the implicit equation of the given
curve is
$$
2-3y-x+2xy=0.
$$

\bigskip
\bigskip

\section{The use of resultants and the Vandermonde matrix}

The computation of the implicit equation of a rationally
parametrized curve can be carried out by computing the {\sl
resultant} of two polynomials. The following theorem provides the
way of doing it (see [15] and [21]).

\medskip
{\bf Theorem 2.} {\it Let
$P=\Big(x(t)={{u_1(t)}\over{v_1(t)}},y(t)={{u_2(t)}\over{v_2(t)}}\Big)$
be a proper rational parametrization of an irreducible curve $C$,
with $gcd(u_1,v_1) = gcd(u_2,v_2) = 1$. Then the polynomial
defining $C$ is  $Res_t(u_1(t) - xv_1(t), u_2(t) - yv_2(t))$ (the
resultant with respect to t of the polynomials $u_1(t) - xv_1(t)$
and $u_2(t) - yv_2(t)$).}
\medskip

Taking this theorem into account, the polynomial defining the
implicit equation of the curve of the example introduced in
Section 1 is
$$
Res_t((1+t)-x(2+t), (3+t)-y(4+t)),
$$
which is precisely
$$
2-3y-x+2xy.
$$
The resultant of two polynomials can be computed, for example, by
using the command {\tt resultant} of the symbolic computation
system {\it Maple}.

However, let us point out here that, in general, the computation
of $Res_t(u_1(t) - xv_1(t), u_2(t) - yv_2(t))$ is not a trivial
task. It is the determinant of the Sylvester (or B\'ezout) matrix
of $u_1(t) - xv_1(t)$ and $u_2(t) - yv_2(t)$ (see [20], for
example), and therefore its computation requires the expansion of
the determinant of a matrix whose entries are polynomials in the
variables $x$ and $y$. The symbolic expansion of a determinant is
a computer algebra problem which requires a lot of time and space
to be solved due to the problem of {\it intermediate expression
swell}. This fact is recognized in [13] where it is said that the
{\it bottleneck} of the algorithm for implicitizing rational
surfaces is the symbolic expansion of the determinant. As it can
be read in [5], one of the most interesting approaches for the
symbolic expansion of the determinant is based on interpolation,
and it is the one presented in [13] and [14].

In that approach the $N$ interpolation data are obtained by
evaluating the Sylvester (or B\'ezout) matrix at the $N$
interpolation nodes
$$
\{(p_1^k, p_2^k): \quad k=0, \ldots, N-1; \quad p_1, p_2 \quad
\text{distinct primes}  \},
$$
and computing the corresponding constant determinants. This clever
selection of the nodes reduces the solution of the interpolation
problem to the solution of a linear system of order $N$ whose
coefficient matrix is the transpose of a nonsingular (since
different monomials evaluate to different values) Vandermonde
matrix (see, for example, [17] for the expression of the
Vandermonde matrix).

\medskip
We illustrate this approach with our example.

We start by showing the Sylvester matrix of the polynomials
$p(t)=1+t-x(2+t)$ and $q(t)=3+t-y(4+t)$, which is
$$S =
\begin{pmatrix}
1-x & 1-2x\\
1-y & 3-4y
\end{pmatrix}.
$$
In order to compute its determinant by using this interpolation
approach, we consider $p_1=2$ and $p_2=3$, and so the $N=4$
interpolation nodes are:
$$
\{ (1,1), (2,3), (4, 9), (8, 27) \}.
$$
The vector with the interpolation data in the corresponding order
is
$$
b=(0, 3, 43, 345)^T,
$$
and the coefficient matrix of the linear system is
$$
A =
\begin{pmatrix}
1 & 1 & 1 & 1 \\
1 & 3 & 2 & 6 \\
1 & 9 & 4 & 36 \\
1 & 27 & 8 & 216
\end{pmatrix}.
$$
The solution of the linear system $Ac=b$ is the vector
$$
c=(2, -3, -1, 2)^T,
$$
and therefore the implicit equation of the curve is
$$
2-3y-x+2xy=0.
$$

In [13], where this approach is presented for the case of surface
implicitization, it is indicated that in this way the problem
reduces to interpolating a univariate polynomial. However, it must
be observed that the coefficient matrix of the linear system is
not a Vandermonde matrix associated with a univariate Lagrange
interpolation problem: it is the transpose of such a matrix.

The approach of [13] has its roots in computer algebra (see, for
example, [23]), and so the solution of Vandermonde systems by
Bj\"orck-Pereyra algorithms [1] is not considered there. A recent
extension of the Bj\"orck-Pereyra approach (closely related to
bidiagonal factorizations of the inverse of the matrix) to
Vandermonde-like matrices can be seen in Chapter 22 of [9].

\bigskip
\bigskip

\section{A choice of nodes leading to the Kronecker product}

Another approach for computing the implicit equation of a plane
rationally parametrized curve by means of resultants and
interpolation is the one introduced in [15].

In that paper, the polynomial defining the implicit equation of
the curve, that is, $Res_t(u_1(t) - xv_1(t), u_2(t) - yv_2(t))$,
is computed by using a bivariate interpolation technique in which
the nodes are arranged forming a tensor product grid. This choice
of the interpolation nodes is specially appropriate for the
interpolation space we are working with, because it reduces the
solution of the interpolation problem to the solution of a linear
system of order $N$ whose coefficient matrix is the Kronecker
product
$$
A = V_x \otimes V_y, $$ (where the Kronecker product $B \otimes D$
is defined by blocks as $(b_{kl}D)$, with $B=(b_{kl})$) with $V_x$
being the Vandermonde matrix  generated by the first component of
the interpolation nodes and $V_y$ being the Vandermonde matrix
generated by the second component of the interpolation nodes (see
[15] for the details). In addition, the algorithm included there
reduces the solution of this linear system with Kronecker product
structure to solving $m+1$ Vandermonde linear systems with the
same matrix $V_y$ and $n+1$ Vandermonde linear systems with the
same matrix $V_x$. In this way, the solution of a bivariate
interpolation problem is reduced to the solution of only
univariate interpolation problems in the variables $x$ and $y$.

A {\it Maple} implementation of the complete algorithm is also
included in the paper. As every linear system to be solved is a
Vandermonde linear system, the algorithm uses the Bj\"orck-Pereyra
algorithm [1, 8] for solving them, since it takes advantage of the
special structure of the coefficient matrices $V_x$ and $V_y$.

The specific choice of the interpolation nodes proposed in the
paper is
$$
\{ (x_i, y_j) = (i, j): i= 0, \ldots, m; j=0; \ldots, n \},
$$
in the same order as the interpolation basis.

\medskip
Now we apply this technique to our example.

\noindent The interpolation nodes considered in this case are:
$$
\{ (0,0), (0,1), (1,0), (1,1) \}.
$$
The vector with the interpolation data in the corresponding order
is:
$$
b=( 2, -1, 1, 0)^T.
$$
The coefficient matrix of the linear system is:
$$A=
\begin{pmatrix}
1 & 0 & 0 & 0 \\
1 & 1 & 0 & 0 \\
1 & 0 & 1 & 0 \\
1 & 1 & 1 & 1
\end{pmatrix}=
\begin{pmatrix}
1 & 0 \\
1 & 1
\end{pmatrix}
\otimes
\begin{pmatrix}
1 & 0 \\
1 & 1
\end{pmatrix}.
$$
The solution of the linear system $Ac=b$ is the vector
$$
c=(2, -3, -1, 2)^T,
$$
and therefore, the implicit equation of the curve is
$$
2-3y-x+2xy=0.
$$

\bigskip
\bigskip

\section{Some remarks on the computational complexity}

We start this section by presenting a bigger (but still small)
example of curve implicitization in which we will show the
different behaviour of the three different approaches we have
described in the previous sections.

\medskip
Let
$$
P=(x(t), y(t)) = \Big(\frac{2t^2+2t+1}{t^3+5},
\frac{t^3-3t^2+t-1}{t^2-3}\Big)
$$
be a rational parametrization of a curve $C$ whose implicit
equation is given by the polynomial
$$
F(x,y)=-53+42y-74y^2+172x+ 707xy+121xy^2+37xy^3-652x^2-1156x^2y
$$
$$
-490x^2y^2-34x^2y^3 +626x^3+396x^3y+432x^3y^2-2x^3y^3.
$$

\medskip
When computing $F(x,y)$ by means of the approach described in
Section $2$, the implicitization problem is reduced to the
solution of a homogeneous linear system of order $16$ with a
non-structured coefficient matrix $A$.  As this matrix is too big
for including it here, we just show one of its largest entries:
$$
A_{15,16}=\frac{761421163154846949}{149346877368718693}.
$$
After solving this linear system the vector with the coefficients
of $F(x,y)$ is obtained.

\medskip
The other two methods need the Sylvester matrix of the polynomials
$p(t)=2t^2+2t+1-x(t^3+5)$ and $q(t)=t^3-3t^2+t-1-y(t^2-3)$ for
computing the implicit equation of $C$. This Sylvester matrix is:
$$
S=
\begin{pmatrix}
-x & 2 & 2 & -5x+1 & 0 & 0\\
0 & -x & 2 & 2 & -5x+1 & 0\\
0 & 0 & -x & 2 & 2 & -5x+1\\
1 & -y-3 & 1 & 3y-1 & 0 & 0\\
0 & 1 & -y-3 & 1 & 3y-1 & 0\\
0 & 0 & 1 & -y-3 & 1 & 3y-1
\end{pmatrix}.
$$

When using the approach described in Section $3$, the
implicitization problem is reduced to the solution of a
non-homogeneous linear system $Ac=b$ of order 16 where $A$ is the
transpose of a Vandermonde matrix. As we have done before, we only
present the greatest entry of the matrix, which corresponds to the
evaluation of the monomial $x^3y^3$ at the interpolation node
$(2^{15}, 3^{15})$:
$$
A_{16,16}=103945637534048876111514866313854976.
$$
As for the vector $b$ with the interpolation data, its largest
component is
$$
b_{16}=-207995995871362988895940143529893921.
$$
The solution of this linear system in which very large numbers are
involved is the vector with the coefficients of $F(x,y)$.

\medskip
Finally, when we compute the implicit equation of $C$ by means of
the method described in Section $4$, the problem is reduced to the
solution of the non-homogeneous linear system $Ac=b$ of order 16
where $A$ is the Kronecker product of two Vandermonde matrices of
order $4$,
$$
A=
\begin{pmatrix}
1 & 0 & 0 & 0\\
1 & 1 & 1 & 1\\
1 & 2 & 4 & 8\\
1 & 3 & 9 & 27
\end{pmatrix}
\otimes
\begin{pmatrix}
1 & 0 & 0 & 0\\
1 & 1 & 1 & 1\\
1 & 2 & 4 & 8\\
1 & 3 & 9 & 27
\end{pmatrix},
$$
and the vector $b$ containing the interpolation data is
{\scriptsize
$$(-53, -85, -265, -593, 93, 72, 35, -12, 2691, 4277, 8723,
15561, 11497, 21242, 44579, 80014)^T.
$$}
In this case only $8$ linear systems with the Vandermonde matrix
of order $4$ presented above have to be solved for obtaining the
coefficients of $F(x,y)$.

\bigskip {\bf Remark.} Let us observe here that, as it can be read
in [20], the implicit equation of a rational curve can also be
computed by using techniques based on computing Gr\"obner bases
with the pure lexicographical ordering. However, although these
techniques are very important from a theoretical point of view,
they are not so effective in practice because they are very time
and space consuming, even for problems of moderate degree. An
example of this situation is the example introduced in this
section. When we tried to compute the implicit equation of $C$ by
using the {\it Maple} command for computing Gr\"obner bases {\tt
gbasis} no answer was obtained after a lot of minutes of
computation, and in addition a lot of memory space (more that 70
megabytes without finishing the computation) was required.

In this sense, it can be read in [7] that the complexity theory
for Gr\"obner bases {\sl gives rise to the pessimistic view that
these methods for polynomial ideals are not useful in practice,
except for rather small cases}. In [20] it is recognized that the
use of Gr\"obner bases for surface implicitization {\sl is not
very computationally efficient}, and this fact is also observed in
a different application of Gr\"obner bases in [3]. As illustrated
with our example in this section, the problems with Gr\"obner
bases appear also in the case of curve implicitization and with
small degrees.

\bigskip As for the approaches based on interpolation, it must be
observed that the numbers involved in the first and in the second
approach are much larger than the numbers involved in the third
one (see the example at the beginning of this section), which can
make the computations slower.

\medskip
Now we briefly analyze the computational complexity of the three
methods for curve implicitization presented in this paper,
beginning with the stage corresponding to the solution of the
linear systems. We will do it in terms of arithmetic operations,
and for the sake of simplicity we assume $m=n$.

\begin{itemize}

\item[-] The computational complexity of the first approach is the
computational complexity of solving a non-structured linear system
of order $N$: $O(N^3)$.

\item[-] The computational complexity of the second approach
is the computational complexity of solving a linear system whose
coefficient matrix is the transpose of a Vandermonde matrix of
order $N$: $O(N^2)$ if the Bj\"orck-Pereyra algorithm is used [1,
8].

\item[-] The computational complexity of the third approach is the
computational complexity of solving $2(n+1)$ Vandermonde linear
systems of order $(n+1)$:  $O((n+1)^3)=O(N^{3/2})$ if the
Bj\"orck-Pereyra algorithm is used for solving each Vandermonde
linear system [1, 8].

\end{itemize}

Let us point out that in the first approach, in addition to the
higher computational cost of solving the linear system, first it
is necessary to evaluate the parametric equations to obtain the
interpolation nodes, and then the coefficient matrix of the linear
system (which is not a structured matrix) must be constructed and
stored.

As for the computation of the interpolation data, its complexity
is the same both in the second and in the third approach, and it
is the complexity of computing $(n+1)^2$ determinants of order
$O(n)$: $O(n^5)=O(N^{5/2})$. Let us point out here that, although
this complexity is greater than the complexity of solving the
linear systems, this stage can easily be parallelized because each
datum can be computed separately.

\medskip
The design of {\sl parallel algorithms} for computer algebra
problems such as resultant computation and elimination of
variables has been recently considered in [2] and [10]. In this
sense, let us observe that the approach introduced in Section 4
has a high degree of intrinsic parallelism present not only in the
computation of the interpolation data but also in the computation
of the coefficients of the interpolating polynomial, since each
one of the linear systems with matrix $V_y$ can be solved
independently and the same happens to each one of the linear
systems with matrix $V_x$ (see [15]).

\medskip
Although it is true that this measure of computational complexity
-which is the standard one in numerical computations- is not a
complete measure of the computational complexity in exact
arithmetic, where the computational cost also depends on the size
and structure of the numbers, it nevertheless can serve to have
useful estimates of the computational complexity. In this sense,
we can read in [19] how the reduction of the order of the
resultant matrix from $2n$ to $n$ may lead to faster computations.

In addition, it must be pointed out that the two implicitization
methods based on interpolation and resultants we have described
can also take advantage of this reduction in the order of the
resultant matrix because they can be used with any other resultant
matrix, for example, the resultant matrix obtained when using the
method of moving curves [19].

\medskip
Finally, it must be noticed that in the situation described in the
approach of Section 4, that is, when the nodes are arranged
forming a Cartesian product grid, explicit formulas (involving the
{\sl Lagrange basis}) exist for computing the interpolation
polynomial [18].

Nevertheless, the existence of an explicit formula does not imply
there is no computational cost in applying it. In this sense we
can recall the classical paper [6], where it is stressed that the
fact that we have Cramer's rule does not make the practical
solution of linear systems a {\sl trivial and dull} task.

In our case, it must be taken into account that the computation of
the interpolation polynomial {\sl in the monomial basis} by using
an explicit formula has a computational cost (see [7] for the cost
in the univariate case). It must be observed that our approach
described in Section 4 has a complexity of the same order as the
algorithm presented in a more general setting in [16].

\bigskip
\bigskip

{\bf Acknowledgements.} This research has been partially supported
by the Spanish Research Grant BFM 2003-03510.

\end{document}